\documentclass{amsart}
\usepackage{amssymb,amsmath,eepic,verbatim}

\newtheorem{theorem}{Theorem}[section]
\newtheorem{corollary}[theorem]{Corollary}

\newtheorem{proposition}[theorem]{Proposition}
\newtheorem{remark}[theorem]{Remark}

\begin{document}
\title{Adjoints and Max Noether's Fundamentalsatz}
\author{William Fulton}
\thanks{Research partially supported by 
NSF Grant DMS9970435}
\address{University of Michigan, Ann Arbor, MI 04109}
\email{wfulton@umich.edu} 
\dedicatory{For S. S. Abhyankar on his 70th birthday}

\begin{abstract}
We give an exposition of the theory of adjoints and conductors for curves 
on nonsingular surfaces, emphasizing the case of plane curves, for which 
the presentation is particularly elementary.  This is closely related to Max 
Noether's ``$AF+BG$'' theorem, which is proved for curves with arbitrary 
multiple components.
\end{abstract}

\maketitle

\section*{Introduction}
Our purpose here is to give an elementary exposition of the theory of adjoints of 
curves in the plane or on a nonsingular algebraic surface.  The treatments we 
have found in the literature are either computationally difficult (\cite{G}, 
\cite{Sa}), or involve quite a bit of advanced machinery: at least the machinery 
of sheaves and cohomology (\cite{Sz}), or even residues and duality (\cite{LJ}).  
See Serre \cite{Se}, Chap. IV, \S12 for a discussion, and Kunz \cite{K} for a 
self-contained treatment of duality in this context.  In addition, we have not 
found Max Noether's ``$AF+BG$'' theorem in its natural generality, which 
allows the given curves to have irreducible components with arbitrary 
multiplicities, so we have taken this opportunity to supply a statement and 
proof.
 	
In honor of Abhyankar, we have made it our goal to be explicit and elementary.  
We have attempted to make this understandable to one who knows only the 
basics of plane curves, and is equipped with an undergraduate algebra course, 
as in \cite{F}, for example; indeed, we expect to include a version of this 
exposition in a revision of \cite{F}.  The case of positive characteristic is included.  The 
local theory applies equally to curves on any nonsingular surface, 
using the less elementary fact 
that the local ring of a point on a nonsingular surface is a unique 
factorization domain.

We thank Joe Lipman for stimulating advice.

\section{Basic concepts and notation}\label{S1}

We work over a fixed algebraically closed ground field  $k$.  We are concerned 
with a nonsingular surface $U$,  a curve $C$ on $U$,  and a point $P$ on $C$.   
We do not assume $C$ is irreducible or reduced, so it may have several 
irreducible components, each appearing with arbitrary positive 
multiplicity.  For a local study, we are free to replace $U$  by any smaller 
neighborhood of $P$,   which we often do without changing notation.  For 
example, we may assume $U$  is affine, and $C$ is defined by an element $F$ 
of the coordinate ring $\Gamma$  of $U$; $F$ is determined up to 
multiplication by a unit.  Instead of working on $U$  and with the coordinate 
rings $\Gamma$ and $\Gamma/(F)$  of $U$  and $C$, we usually work with 
the limiting local rings  $\Lambda = \mathcal{O}_PU$,  the local ring of $U$  
at $P$, and  $A = \Lambda/(F) = \mathcal{O}_PC$,  the local ring of $C$ at 
$P$. Let  $\mathfrak{M}$  denote the maximal ideal of  $\Lambda$,  and  
$\mathfrak{m}$  the maximal ideal of  $A$.  
	
The {\em multiplicity} of $C$ at $P$ is the largest integer  $r = r_P(C)$  such 
that $F$ is in  $\mathfrak{M}^r$.   If $x$ and $y$ generate  $\mathfrak{M}$,  
the image  $F_r$  of $F$ in  $\mathfrak{M}^r/\mathfrak{M}^{r+1}$  can be 
written  $F_r(x,y) = \sum_{i+j = r} a_{ij}x^iy^j$.  The irreducible factors of  
$F_r$  give the {\em tangent lines} to $C$ or 
$F$ at $P$. We say that the coordinates 
$x$ and $y$ are {\em suitable} for 
$C$ or $F$ if  $F_r(0,1) \neq 0$.  This can always 
be achieved by a linear change of coordinates.  In this case  $F_r = a\prod(y - 
\alpha_ix)^{m(i)}$,  with $a \neq 0$,  and distinct  $\alpha_i$  in  $k$; the 
tangent lines are the lines  $y = \alpha_ix$.   
	
For simplicity, we start with (and readers who wish may remain with) the case 
where $C$ is planar, by which we mean that, after shrinking $U$  if necessary, 
$U$  is isomorphic to an open subset of the plane  $\mathbb{A}^2$,  and $C$ 
is defined by a polynomial $F = F(x,y)$  in  $k[x,y]$.  Applying a translation, we 
may assume $P$ corresponds to the origin  $(0,0)$.  Then $F =  F_r + F_{r+1} 
+ \ldots+ F_n$  where each $F_d$  is a homogeneous polynomial of degree  
$d$  in $x$ and $y$, and  $F_r \neq 0$.
	
The blow-up of $U$  at $P$ can be described as follows.  Shrinking $U$  if 
necessary, we may assume there are $x$ and $y$ in $\Gamma$ that generate 
the maximal ideal of $P$. The {\em blow-up} $U'$  of $U$  at $P$ is the 
subvariety of $U \times \mathbb{P}^1$  defined by the equation  $xT = yS$,  
where  $S$  and  $T$  are the homogeneous coordinates on  $\mathbb{P}^1$.  
The {\em exceptional divisor}  $E \cong \mathbb{P}^1$ is defined by  $x = y = 
0$.  The blow-up is covered by two affine open subsets  $U'_0$  and  $U'_1$,  
where  $S$  and  $T$,  respectively, are not zero.  The first,  $U'_0$,  is the 
subvariety of $U\times \mathbb{A}^1$  defined by  
$y = tx$,  where $t$  ($= T/S$)  is the coordinate on  $\mathbb{A}^1$;  
similarly,  $U'_1$ is the subvariety of $U\times \mathbb{A}^1$  defined by  $x 
= sy$,  where  $s$  ($= S/T$)  is the coordinate on  $\mathbb{A}^1$.  Note 
that if $U$  is planar, both  $U'_0$  and  $U'_1$ are planar; indeed, when 
$U= \mathbb{A}^2$,  $U'_0$  and  $U'_1$  are both isomorphic to  
$\mathbb{A}^2$,  by the maps  $(x,y,t) \mapsto (x,t)$  and  $(x,y,s) \mapsto 
(y,s)$.  
	
The exceptional divisor  $E$  is defined on  $U'_0$  by  $x$,  and on  $U'_1$  
by  $y$.  The coordinate functions $x$ and $y$ determine a basis for  
$\mathfrak{M}/\mathfrak{M}^2$,  where $\mathfrak{M}$  is the maximal 
ideal of $\Gamma$;  this determines an identification of  $E$  with the 
projective tangent 
space to $U$ at $P$. (This can also be seen by the intrinsic construction of the 
blow-up as  $\operatorname{Proj}(\oplus \mathfrak{M}^i)$,  with  $E = 
\operatorname{Proj}(\oplus \mathfrak{M}^i/\mathfrak{M}^{i+1})$,  but we do 
not need this description.)
	
The {\em proper transform}  $C'$  of $C$ is the curve on  $U'$  defined by the 
equation  $\pi^*(C) = C' + r E$,  where  $r$  is the multiplicity of $C$ at $P$.  
Explicitly,  $C'$  is defined on  $U'_0$  by the $F'$  such that $F = x^rF'$,  
and  $C'$  is defined on  $U'_1$   by an  $F''$  with $F = y^rF''$.  When $U$ is 
planar, on  $U'_0$ we have   $F(x,y) = F(x,xt) = x^rF'$,  where 
\[
F'  =  F_r(1,t) + xF_{r+1}(1,t) + \ldots+ x^{n-r}F_n(1,t)\,,
\]
and similarly on  $U'_1$.  
	
The coordinates $x$ and $y$ are  suitable for $C$ if the line 
$x = 0$  is not tangent to $C$ at $P$. This means that  $C'$  does not contain 
the point  $[0:1]$  in  $\mathbb{P}^1 = E$,  i.e., $ C'$  is contained in the 
affine piece  $U'_0$.  In this case we set  $\Lambda' = \Lambda[t]/(y-xt)$,  and  
$A' = \Lambda'/(F') = \Lambda[t]/(y-tx, F')$.  Note that  
$\mathfrak{M}\Lambda' = (x)$.  The projection from $C'$  to $C$ corresponds 
to the natural homomorphism from $A = \Lambda/(F)$  to  $A'$.  We often 
abuse notation by writing  $x$, $y$, and  $t$  for their images in $A$  or  $A'$.  
(In fancier and more intrinsic language, 
the morphism  $\pi \colon C' \to C$  is a 
finite morphism, and  $A'$ is the 
localization of  $\pi_*(\mathcal{O}_{C'})$  at 
$P$.)  We will always assume that coordinates 
are suitable for any finite number 
of given curves passing through  $P$.
	
We use the fact that  $\Lambda$  is a unique factorization domain, which is a 
general fact about regular local rings (cf. \cite{M1}, \S19 or \cite{E}, \S19.4); 
in the planar case it follows from the fact that the polynomial ring  $k[x,y]$  is a 
unique factorization domain.  If $G$  is an element of  $\Lambda$  that has no 
irreducible factors in common with $F$, it follows that the image of $G$  in  $A 
= \Lambda/(F)$  is a non-zero-divisor.  The same holds for local rings on the 
blow-up, and for  $\Lambda'$.  In particular, if two elements $G$  and $H$ of  
$\Lambda'$  have only a finite number of common zeros, then each is a 
non-zero 
divisor in the ring modulo the ideal generated by the other.  For example, 
$x$ is a non-zero-divisor in $A$  and in  $A'$;  in the latter case this follows 
from the fact that the common zeros of $x$ and $F'$ correspond to the finite 
number of tangent lines to $C$ at $P$.

\section{Conductors}\label{S2}

Recall that for any subring $A$  of  a ring  $A'$,  the {\em conductor} $I$ of 
$A$  in  $A'$ is the ideal of elements  $a$  in $A$  such that  $aA' \subset A$;  
it is the largest ideal of $A$  that is also an ideal in  $A'$;  and any element  $a$  
in $I$ satisfies  $aA' \subset I$.
	
Our results depend on the following elementary computation.    

\begin{proposition}\label{P2.1} Suppose  $A \to A'$ arises as in Section 1 from 
the blowup of a curve at a point.  
\begin{enumerate}
\item[(i)] The homomorphism  $A \to A'$ is injective;  $A'$ is a finitely 
generated  $A$-module, generated by the elements  $t^j$,  $0 \leq j\leq  r-1$.  
\item[(ii)]  The images of the elements  $x^it^j$,  for  
$0 \leq i < j \leq r-1$,  form a 
basis for  $A'/A$  over  $k$.  
\item[(iii)] The conductor $I$ of $A$  in  $A'$ is  $\mathfrak{m}^{r-1} = x^{r-
1}A'$;  the images of the elements  $x^iy^j$,  for  $0 \leq i + j\leq  r-2$,  
form a 
basis for $A/I$ over  $k$.
\end{enumerate}
\end{proposition}
\begin{proof}  For any $H$ in  $\Lambda'$,  there is a positive  $N$  such that  
$x^NH$  is in  $\Lambda$.  It follows that if $G$  in  $\Lambda$  is not 
divisible by $F$, then the image of $G$  in  $\Lambda'$  is not divisible by  
$F'$;  for if  $G = F'H$  in  $\Lambda'$,  then  $x^NG = FJ$  in  $\Lambda$  
for  some integer  $N$  and some element  $J$  in  $A$;  since $x$ is not a 
zero-divisor in  $A$,  this is a contradiction.  This shows that $A$  is a subring 
of  $A'$.  To show that  $A'$ is generated over $A$  by the  $r$  elements  $1$, 
$t$, $t^2$, \ldots, $t^{r-1}$,  by Nakayama's Lemma (cf. \cite{E}, \S4.1 or 
\cite{N}, \S{}I.4), 
it suffices to show that  $A'/xA'$  is generated over  $A/xA$  
by the images of these elements.  But  $A'/xA' = \Lambda[t]/(x,y,F'_r(1,t)) = 
k[t]/(F'_r(1,t))$,  and  $F'_r(1,t)$  is a polynomial in  $t$  of degree  $r$.  This 
proves (i).  
	
Note that $x^{r-1}t^j = x^{r-1-j}y^j$  for  $j \leq r-1$,  so  $x^{r-1}$  is in the 
conductor $I$. And since  $\mathfrak{m} = (x,y) \subset xA'$,  we have  
$\mathfrak{m}^{r-1} \subset x^{r-1}A' \subset I$.  The other assertions come 
from looking at the exact sequence  
\[
A/\mathfrak{m}^{r-1}  \to  A'/x^{r-1}A'  \to  A'/A  \to  0\, ,
\]
where the second map is the canonical surjection, arising from the fact that  
$x^{r-1}A' \subset A$;  and the first is induced from the inclusion of $A$  in  
$A'$,  noting that  $\mathfrak{m}^{r-1} \subset x^{r-1}A'$.   Since $F$ is in  
$\mathfrak{m}^r$,  we know that the images of the elements  $x^iy^j$,  for  $i 
+ j < r-1$,  form a basis for  $A/\mathfrak{m}^{r-1} \cong 
\Lambda/\mathfrak{m}^{r-1}$.   We claim that the images of the elements  
$x^it^j$,  for  $i < r-1$  and  $j < r$,  form a basis for  $A'/x^{r-1}A'$.  Since 
$x$ is a non-zero-divisor in  $A'$,  looking at the filtration  $A' \supset xA' 
\supset x^2A' \supset \ldots \supset x^{r-1}A'$,  it suffices to show that the 
elements  $t^j$,  $j < r$,  form a basis for  $A'/xA'$.  And this is clear since, as 
we have seen,  $A'/xA' = k[t]/(F'(1,t))$.  
	
The mapping from  $A/\mathfrak{m}^{r-1}$  to  $A'/x^{r-1}A'$  takes  
$x^iy^j$  to  $x^{i+j}t^j$.   It follows that the images of the remaining  
$x^it^j$,  with  $i < j < r$,  form a basis for  $A'/A$, which proves (ii).  (We also 
see that the displayed sequence is exact on the left.)

To finish the proof of (iii), we must show that  $I \subset \mathfrak{m}^{r-1}$.  
Since we know that  $\mathfrak{m}^{r-1} \subset I$,  if this were not true 
there would be an element  $z$  in $I$ of the form  $\sum_{i+j < r-1} 
a_{ij}x^iy^j$,  with each  $a_{ij}$  in  $k$,  and not all  $a_{ij} = 0$.  Let  $\ell$  
be minimal such that some  $a_{\ell j} \neq 0$.  Then 
\begin{align}
z \cdot t^{\ell+1}  &=  \sum_j a_{\ell j}x^{\ell}(xt)^jt^{\ell+1}  
+  \sum_{i > 
\ell} \sum_j a_{ij}x^iy^jt^{\ell+1} \notag \\
{}&=  \sum_j a_{\ell j}x^{\ell+j}t^{\ell+j+1}  
+  \sum_{i > 
\ell} \sum_j a_{ij}x^{i-\ell-1}y^{j+\ell+1} \,. \notag   
\end{align}
Since  $z\cdot t^{\ell+1}$  is in  $A$,  and the second term on the right is in  
$A$,  the first term must also be in  $A$.  But since each  $\ell+j+1 < r$,  it 
follows from (ii) that no such linear combination can be in  $A$. 
\end{proof}
  
\begin{corollary}\label{C2.2} The dimensions (over  $k$)  of  $A'/A$  and  
$A/I$  are both equal to  $r(r-1)/2$.  
\end{corollary}

\begin{corollary}\label{C2.3}  The image of any non-zero-divisor in $A$  is a 
non-zero-divisor in  $A'$. 
\end{corollary}

\begin{proof}  Suppose  $a$  is in  $A$,  and  $a\cdot a' = 0$  for some  $a'$  
in  $A'$.  Then  $a\cdot x^{r-1}\cdot a' = 0$,  and  $x^{r-1}\cdot a'$  is in  
$A$.  If  $a$  is a non-zero-divisor, then  $x^{r-1}\cdot a' = 0$.  But we have 
seen that $x$ is a non-zero-divisor in  $A'$,  so  $a' = 0$. 
\end{proof}

These results extend readily to the case of the blow-up of several points.  This 
case is not quite local, so one needs a slightly more general blowing up.  If  $V$  
is a nonsingular surface, embedded as a locally closed subvariety in some 
projective space  $\mathbb{P}^n$,  and $P$ is a point in  $V$,  one may 
choose 
homogeneous linear polynomials 
 $L_0, \ldots, L_n$  whose restrictions to $V$ vanish only 
at $P$. Then  $[L_0 : \ldots : L_n]$  determines a morphism from  $V 
\!\smallsetminus\! {P}$  to  $\mathbb{P}^n$.  The blowup  $V'$  of  $V$  at 
$P$ can be taken to be the closure of the graph of this morphism in  $V \times 
\mathbb{P}^n$.  (This shows that  $V'$  can also be embedded as a locally 
closed subvariety of a projective space, since by Segre  
$\mathbb{P}^n \times \mathbb{P}^n$  is a closed subvariety of a larger 
projective space.)  The projection from  $V'$  to  $V$  is an isomorphism over  
$V \!\smallsetminus\! {P}$.  To see that it is isomorphic to the blowup 
considered before, over some affine neighborhood $U$ of $P$,  take such a 
neighborhood, with functions $x$ and $y$ generating the maximal ideal of $P$. 
Let  $L$  be a linear form that does not vanish at $P$,   and write  $L_i/L = 
a_{i1}x + a_{i2}y$  for some functions  $a_{ij}$.  Shrink $U$ if necessary so 
that the matrix  $(a_{ij})$  has rank  $2$  everywhere on $U$.  Then this 
matrix  $(a_{ij})$  determines a closed embedding of $U \times 
\mathbb{P}^1$  in $U \times \mathbb{P}^n$,  and one verifies easily that the 
blowup we defined earlier in $U \times \mathbb{P}^1$  is mapped to the 
closure of the graph just defined.  
	
Now suppose $U$ is an affine nonsingular surface, and  $P_1, \ldots, P_s$  are 
distinct points of $U$. In this case we take  $\Lambda$  to be the semi-local 
ring that is the localization of the coordinate ring $\Gamma$ of $U$ at the 
multiplicative set of elements not vanishing at any  $P_i$.  If $C$ is a curve on 
$U$ (with irreducible components of arbitrary multiplicities), set  $A = 
\Lambda/I(C)$,  where  $I(C)$  is the ideal of elements which are divisible by a 
local equation for $C$ at each  $P_i$.  Let  $U' \to U$  be the simultaneous 
blow-up of $U$ at each of the points  $P_i$  (i.e., the result of successively 
blowing up each  $P_i$,  the result being independent of the order of blow-up).  
We again have the proper transform  $C'$  of $C$,   with its finite morphism  
$C' \to C$,  which corresponds to a monomorphism $A \to A'$ of  
$k$-algebras.  If $I$ is the conductor, then  $A/I$  and $A'/A$ both have 
dimension  $\sum_i r_i(r_i-1)/2$,  where  $r_i$  is the multiplicity of $C$ at  
$P_i$.  The point is that, since the  $A$-modules and $A'/A$ have support at 
these points  $P_i$,  we have canonical decompositions  $A/I = \oplus_i 
A_i/IA_i$,  and  $A'/A = \oplus_i A'_i/A_i$,  where  $A_i = 
S_i^{-1}A$,  with  $S_i$  the multiplicative set of elements in $A$  not 
vanishing at  $P_i$,  and  $A'_i = S_i^{-1}A'$  (cf. \cite{F}, \S2.9, 
\cite{E}, \S2.4).  Each of  $A_i \to A'_i$  
is an extension as studied above, so 
we know the dimensions of each summand in these decompositions.
	
Because of this we may repeat the blowing up process.  Starting from the 
blow-up
 $U^{(1)} = U' \to U^{(0)} = U$  of $U$ at $P$,  one can construct the 
blow-up 
$U^{(2)} \to U^{(1)}$  of  $U^{(1)}$  at a finite number of points mapping 
to $P$ (lying on the exceptional divisor).  Repeating, at each stage blowing up 
points in the exceptional divisors from the preceding stage, we get a sequence 
\[
U^{(n)} \to U^{(n-1)} \to \ldots \to U^{(2)} \to U^{(1)} = U' \to 
U^{(0)} = U \, .
\]
Points in any  $U^{(n)}$  mapping to $P$ are called {\em infinitely near} points 
to $P$,   in the  $n^{\text{th}}$  neighborhood, 
for $n \geq 0$.  If at each stage  $C^{(i)}$  is 
the proper transform of  $C^{(i-1)}$,  we have a sequence 
\[
\widetilde{C} = C^{(n)} \to C^{(n-1)} \to \ldots \to C^{(2)} \to C^{(1)} = 
C' \to C^{(0)} = C \, , 
\]
and corresponding finite extensions of  $k$-algebras: 
\[
A = A^{(0)} \subset A' = A^{(1)}  \subset A^{(2)} \subset \ldots \subset  
A^{(n-1)} \subset  A^{(n)} = \widetilde{A} \, .
\]
If  $Q$  is an infinitely near point, in the neighborhood  $U^{(i)}$,  we let  
\[
r_Q  =  r_Q(C)  =  r_Q(F)  
\]
be the multiplicity of the proper transform $C^{(i)}$  at  $Q$. 

\begin{proposition}\label{P2.4}  Let $I$ be the conductor of $A$  in  $A'$,  $J$  
the conductor of  $A'$ in  $\widetilde{A}$,  and $K$  the conductor of $A$  in  
$\widetilde{A}$.  Then  $K = I\cdot J$.
\end{proposition}
\begin{proof}  From the definition of conductors we have  $I \cdot J \subset 
K$;  we must show that  $K \subset I \cdot J$.  Since forming conductors 
commutes with localization, we may assume $A$  is the local ring of one point 
$P$. Choosing coordinates as above, we have seen that  $I = \mathfrak{m}^{r-
1} = x^{r-1}A'$.  If  $u$  is an element of  $K$,  since  $K \subset I$  from the 
definition, we may write  $u = x^{r-1}\cdot v$,  for some  $v$  in  $A'$.  It 
suffices to show that  $v$  is in $J$, i.e., that  $v\cdot b$  is in  $A'$ for any  
$b$  in  $\widetilde{A}$.  Since  $u$  is in  $K$,  $u\cdot b$  is in  $K \subset 
I = x^{r-1}A'$,  so we can write  $u\cdot b = x^{r-1}\cdot a'$  for some  $a'$  
in  $A'$.  But then  $x^{r-1}\cdot v \cdot b = x^{r-1}\cdot a'$.  
By Corollary 
\ref{C2.3}, $x$ is a non-zero-divisor in  $\widetilde{A}$,  and 
it follows that  
$v\cdot b = a'$,  as desired.
\end{proof}

\begin{corollary}\label{C2.5} $\dim(\widetilde{A}/A) = \dim(A/K) = \sum 
r_Q(r_Q-1)$,  the sum over all infinitely near points  $Q$  in some  $U^{(i)}$,  
$0 \leq i \leq n-1$.  
\end{corollary}
\begin{proof}  As in the proposition, we may assume $A$  is local.  We know 
that  $\dim(A'/A) = \dim(A/I) = r_P(r_P-1)/2$.  By induction on the length of 
the chain, we have  $\dim(\widetilde{A}/A') = \dim(A'/J) = \sum r_Q(r_Q-1)$,  
the sum over infinitely near  $Q$  in some  $U^{(i)}$,  
$1 \leq i \leq n-1$.  From the 
inclusions  $A \subset A' \subset \widetilde{A}$  we have  
$\dim(\widetilde{A}/A') = \dim(\widetilde{A}/A') + \dim(A'/A)$.  It therefore 
suffices to show that $\dim(A/K) = \dim(A'/J) + \dim(A/I)$;  adding  
$\dim(A'/A)$  to both sides, we are reduced to proving that  $\dim(A'/K) = 
\dim(A'/J) + \dim(A'/I)$.  
Since $K \subset J$, so $\dim(A'/K) = \dim(A'/J) + \dim(J/K)$, this is 
equivalent to proving that $\dim(J/K) = \dim(A'/I)$.  
Since  $K = x^{r-1}J$  and  $I = x^{r-1}A'$, we 
must show that  $\dim(J/x^{r-1}J) = \dim(A'/x^{r-1}A')$.  From the inclusions 
\[
x^{r-1}J  \subset  J  \subset  A'   \qquad \text{ and } \qquad   
x^{r-1}J  \subset  x^{r-1}A'  \subset  A' \, ,
\]
we have $\dim(J/x^{r-1}J)+\dim(A'/J) = \dim(x^{r-1}A'/x^{r-1}J) 
+\dim(A'/x^{r-1}A')$, so we are reduced to showing that  $\dim(A'/J) = 
\dim(x^{r-1}A'/x^{r-1}J)$.  But multiplication by  $x^{r-1}$  gives an 
isomorphism of  $A'/J$  with  $x^{r-1}A'/x^{r-1}J$, since  $x^{r-1}$  is a 
non-zero-divisor in  $A'$. 
\end{proof}

The fact that  $\dim(\widetilde{A}/A) = \dim(A/K) = 
\frac{1}{2}\dim(\widetilde{A}/A)$  is known as Gorenstein's theorem \cite{G}.  
It depends on the fact that $C$ is a curve on a nonsingular surface.  For 
example, if $C$ is the curve in affine $3$-space which is the image of the map  
$t \mapsto (t^3, t^4, t^5)$  from the affine line, the conductor $K$ of $A = 
k[t^3,\, t^4,\, t^5]$  in $\widetilde{A} = k[t]$  is generated by the maximal ideal 
at the origin, but the images of  $t$  and  $t^2$ give a basis for 
$\widetilde{A}/A$.  The same is true after localizing at the origin, so one has an 
example where  $\dim(A/K) = 1$  but  $\dim(\widetilde{A}/A) = 2$.  For 
another proof that  $\dim(\widetilde{A}/A) = \sum r_Q(r_Q-1)$, see \cite{D}.
 
\begin{corollary}\label{C2.6}  Let $G$  and $H$ be elements of  $\Lambda$,  
with images  $g$ and  $h$  in  $A$.  Let $D$ be the curve defined by $G$, and 
assume that the proper transforms  $D^{(n)}$  and  $C^{(n)}$  in  $U^{(n)}$  
have no points in common.  If  $r_Q(H) \geq r_Q(D) + r_Q(C) - 1$  for all 
infinitely near points  $Q$  to $P$ in $C$,   then  $h$  is in  $g\cdot K$,  with 
$K$  the conductor of $A$  in $\widetilde{A}$.   In particular,  $h$  is divisible 
by  $g$  in  $A$.
\end{corollary}
\begin{proof}  If  $n = 0$,  there is nothing to prove.  Let  $a$  and  $b$  be the 
multiplicities of $G$  and $H$ at $P$. For the first blowup, choosing 
coordinates that are suitable for $G$  and $H$ as well as $F$, we have  
$g = x^{a-1}g'$,  $h = x^{b-1}h'$,  with  $g'$,  $h'$  in  $A'$.  By induction on 
the length of the chain, we know that  $h' = g'\cdot z$,  with  $z$  in $J$. 
Since $b-a-r+1 \geq 0$,
  $x^{b-a-r+1}z$  is in $J$, so  $x^{b-a}z = x^{r-1}(x^{b-a-r+1}z)$  is in  
$x^{r-1}J = K$.  Therefore  
$h = x^{b-1}h' = g'x^{b-1}z = g \cdot x^{b-a}z$  is in  $g \cdot K$,  as required.  
\end{proof}

Now suppose $C$ is an irreducible curve at $P$, so its local ring $A$  is an 
integral domain.  Since $A$  is the localization of a finitely generated algebra 
$\Gamma$ over the field  $k$,  it is a general theorem of E. Noether (see 
\cite{E}, \S13.3 or \cite{N}, \S36) that the integral closure $\widetilde{A}$  of 
$A$  in its quotient field is a finitely generated  $A$-module.  If at each stage of 
blowing up, one blows up at all the singular points in exceptional divisor of each 
proper transform  $C^{(i)}$,  one arrives at a chain  $A = A^{(0)} \subset 
A^{(1)} \subset \ldots \subset A^{(n)} \subset \widetilde{A}$.  It follows that 
this process must terminate, so $\widetilde{A} = A^{(n)}$  for some  $n$.  
Indeed the dimension of $\widetilde{A}/A$  puts a bound on the number of 
steps required.  
	
In the planar case, one can see that this process terminates directly, without 
using Noether's theorem.  We include a proof in the appendix.  
	
Suppose $C$ is irreducible, and one performs the sequence of blowups to 
resolve the singularities of $C$, so  $C^{(n)} = \widetilde{C}$  is nonsingular.  
In this case the conductor $K$  of $A$  in $\widetilde{A}$  is called the {\em 
conductor of} $C$ {\em at} $P$. An element $G$  in  $\Lambda$  is {\em 
adjoint} to $C$ at $P$ if the image of $G$  in $A$  is in the ideal  $K$.  Define 
an effective divisor  $\Delta_P = \sum d_QQ$  on  $\widetilde{C}$  by  
defining 
$d_Q$  to be the order of the ideal $K$  at  $Q$,  i.e.,  $K \cdot 
\mathcal{O}_Q\widetilde{C} = \mathfrak{m}_Q(\widetilde{C})^{d_Q}$.  The 
degree of  $\Delta_P$  is the dimension of $\widetilde{A}/K$,  which, 
by Corollary \ref{C2.5},  is  
$2\cdot \delta_P$,  with  $\delta_P = \dim(\widetilde{A}/A) = \dim(A/K)$.  
An element  $h$  of the function field  $R(C)$,  i.e., the quotient field of  $A$,  
is in the conductor $K$  if and only if  $\operatorname{ord}_Q(h) \geq d_Q$  
for all  $Q$  in  $\widetilde{C}$,  where  $\operatorname{ord}_Q$  is the 
order function on  $R(C)$  defined by the discrete valuation ring  
$\mathcal{O}_Q(\widetilde{C})$.  If  $g$  and  $h$  are in  $R(C)$,  and  
$\operatorname{ord}_Q(h) \geq \operatorname{ord}_Q(g) + d_Q$  for all  $Q$  
in  $\widetilde{C}$,  then  $h$  is in  $g \cdot K \subset g \cdot 
\mathcal{O}_P(C)$.    

\begin{remark}\label{R2.7} For any  $g$  in $A$  
which is a non-zero-divisor in  $A'$, we have   
$$\dim(A'/gA') = \dim(A/gA)$$.  As in Corollary \ref{C2.5}, this follows by 
comparing $gA \subset gA' \subset A'$ and  $gA \subset A \subset A'$,  and 
noting that $A'/A$ is isomorphic to  $gA'/gA$.  In particular, we see that  
$\dim(A'/xA') = \dim(A/xA)$.
\end{remark}
	
This analysis gives a quick proof of the following formula of Max Noether for the 
intersection multiplicity of two curves $C$ and $D$ at $P$. Here we assume 
$C$ and $D$ have no irreducible components in common through $P$. The 
intersection multiplicity  $I(P,C \cdot D)$  is 
defined to be the dimension over  
$k$  of  $\Lambda/(F,G)$,  where $F$ and $G$  are local equations for $C$ and 
$D$.	 

\begin{proposition}\label{P2.8}  The intersection multiplicity is given by the 
formula
\[
I(P,C \cdot D)  =  \sum r_Q(C) \cdot r_Q(D) \,,  
\]
where the sum is over 
$Q = P$ and all infinitely near points  $Q$  of $P$ that lie in proper 
transforms of both $C$ and  $D$.	
\end{proposition}
\begin{proof}  By induction, we need only show that  $I(P,C \cdot D) = 
r_P(C)\cdot r_P(D) + \sum I(P',C'\cdot D')$,  where  P'  varies over the points 
lying over $P$ in both proper transforms  $C'$  and  $D'$,  in the blowup  $U'$  
of $U$ at $P$.  Let  $g$  be the image in  $A = \mathcal{O}_P(C)$  of a local 
equation for $D$ at $P$. By Remark \ref{R2.7},
\[
I(P,C\cdot D)  =  \dim(A/gA)  =  \dim(A'/gA') \, .
\]
In  $A'$,  $g = x^sg'$,  where  $s = r_P(D)$.  There is an exact sequence 
\[
	0  \to  A'/g'A'  \to  A'/x^sg'A'  \to  A'/x^sA'  \to  0 \, ,
\]
where the first map is multiplication by  $x^s$,  and the second is the natural 
projection; the exactness follows from the fact that $x$ is a non-zero-divisor in  
$A'$ (cf. \cite{F}, \S3.3).  Therefore 
\[	
	\dim(A'/gA')  =  \dim(A'/x^sA') + \dim(A'/g'A')\, .
\]
Since  $\dim(A'/x^sA') = \dim(A/x^sA) = s \cdot \dim(A/xA) = s \cdot 
r_P(C)$,  and  $\dim(A'/g'A') = \sum I(P',C'\cdot D')$,  the conclusion follows.
\end{proof}
		 
\begin{corollary}\label{C2.9}  For any $C$ and $D$ with no irreducible 
components in common at $P$,   there is a sequence of blowups so that  
$C^{(n)}$  and  $D^{(n)}$  are disjoint, and each is a disjoint union of curves, 
each consisting of a nonsingular curve with some positive multiplicity.
\end{corollary}
\begin{proof}  We know that a sufficient number of blowups will make each 
irreducible component of $C$ and of $D$ nonsingular.  By the proposition, a 
sufficient number of blowups will then make the proper transform of pairs of 
these components disjoint.
\end{proof}
\begin{remark}\label{R2.10} The same reasoning shows that one can make the 
total transforms, including all the exceptional divisors and their proper 
transforms, a union of nonsingular curves, with some multiplicities, with each 
pair of irreducible components meeting transversally, and no three 
components passing through any point.
\end{remark}

\section{Adjoints and Differentials}\label{S3}

For any $k$-algebra  $R$, we have the $R$-module  $\Omega_{R/k}$  of 
differentials over  $k$.  It can be defined to be the free $R$-module on symbols  
$df$,  for  $f$  in  $R$,  modulo the submodule generated by all:  
1)~$df$, $f$ in $k$;  2)~$d(f+g) - df - dg$,  $f$, $g$  in  $R$;  3)~$d(fg) - 
f\cdot dg - g\cdot df$,  $f$, $g$  in  $R$.   It is constructed so that for any  
$R$-module  $M$,  the $k$-linear derivations from  $R$  to  $M$  correspond 
to $R$-linear homomorphisms from  $\Omega_{R/k}$  to  $M$.   See \cite{F}, 
\S8.4 or \cite{E}, \S16 for basic facts about differentials.
	 
Let $\mathcal{K}$  
be the field of rational functions on $U$,  i.e., the quotient field of  
$\Lambda$.  The differentials  $\Omega_{\mathcal{K}/k}$  
form a vector space over $\mathcal{K}$  
of dimension  $2$.  If $x$ and $y$ generate the maximal ideal of  $\Lambda$,  
then  $dx$  and  $dy$  give a basis for  $\Omega_{\mathcal{K}/k}$  
over  $\mathcal{K}$.  If $F$ is a 
local equation for an irreducible curve $C$,  then  $dF = F_x\,dx + F_y\,dy$  is 
not zero on $U$,  although it vanishes on $C$.   If  $F_y \neq 0$,  then its 
image in $A$  is not zero.  This can be seen by induction on the length of steps 
needed to resolve the singularity, it being clear when $P$ is a nonsingular 
point of $C$.   Starting with the equation $F = x^rF'$,  with  $y = xt$,  
differentiating both sides with respect to  $t$  gives  $xF_y = x^rF'_t$,  i.e.,   
\[
F_y  =  x^{r-1}F'_t \,.  
\]
By induction, we know that the image  $F'_t$  is not zero in  $A'$,  and since 
$x$ is a non-zero-divisor,  $F_y$  is not zero in  $A$. 
	
The differentials  $\Omega_{R(C)/k}$  of the function field of $C$ over  $k$  
form a $1$-dimensional vector space over  $R(C)$.  It is generated by   
$dz$,  where  $z$  is any element of  $R(C)$  such that  $R(C)$  is a finite 
separable extension of  $k(z)$.  We consider the differential 
\[
\omega = dx/F_y  \quad \text{ on }  \quad C \,.   
\]
(If  $F_y = 0$  on $C$,   then we use  $\omega = - dy/F_x$.)  This differential 
on $C$ is independent of choice of coordinates, up to multiplication by a 
function not vanishing at $P$.  Explicitly, if  $x = x(u,v)$,  $y = y(u,v)$,  for  
$u$  and  $v$  other coordinates, and we set  $\widetilde{F}(u,v) = 
F(x(u,v),y(u,v))$,  then, if  $\widetilde{F}_v \neq 0$,  a calculation as in 
calculus shows that 
\[
dx/F_y  =  J \cdot du/\widetilde{F}_v \, ,
\]
where  $J = x_uy_v - x_vy_u$  is the Jacobian determinant; if  
$\widetilde{F}_u \neq 0$,  then  $dx/F_y  =  - J \cdot dv/\widetilde{F}_u$.  
	
Recall that a differential  $\omega$  on a nonsingular curve  $\widetilde{C}$  is 
{\em regular} at a point  $Q$  if, for a uniformizing parameter  $t$  for  
$\widetilde{C}$  at $Q$,  $\omega = h\, dt$,  with  $h$  regular at  $Q$,  i.e.,  
$h$  is in  $\mathcal{O}_Q(\widetilde{C})$.  More generally, the {\em order}  
$\operatorname{ord}_Q(\omega)$  is the order of such a function  $h$.  

\begin{proposition}\label{P3.1}  An element  $g$  in  $A = \mathcal{O}_P(C)$  
is in the conductor $K$  of $\widetilde{A}$  over $A$  if and only if the 
differential  $g \cdot dx/F_y$  is regular at each point of  $\widetilde{C}$  that 
maps to  $P$.
\end{proposition}
\begin{proof}  Take coordinates as before at $P$.  If $P$ is a nonsingular point 
on $C$,  then  $F_y$  or  $F_x$  is a unit at $P$,  and the assertion is clear.  
Otherwise perform a blowup, and write $F = x^rF'$.   As we have just seen,  
$F_y = x^{r-1}F'_t$,  where $x$ and  $t$  are coordinates on  $U'$.  By 
induction on the number of blowups needed to resolve the singularity, a 
function  $g$  is in the conductor $J$  for  $A'$ in $\widetilde{A}$  exactly 
when $g\cdot dx/F'_t$  is regular at all points of  $\widetilde{C}$  over $P$. 
For  $g$  to be in the conductor $K = x^{r-1}J$,  $g/x^{r-1}$  must be in 
$J$, i.e.,  $g\cdot dx/F_y  =  (g/x^{r-1})dx/F'_t$  must be regular at all points 
of  $\widetilde{C}$  over $P$, as required.
\end{proof}

\section{Plane Curves}\label{S4}

The results on adjoints lead to a sharp form of Noether's theorem, allowing 
curves with arbitrary multiple components.

\begin{theorem}[Max Noether's Fundamentalsatz]\label{T}  Let $C$ and $D$ be 
plane curves with no common components, defined by homogeneous 
polynomials $F$ and $G$  of degrees  $c$  and  $d$.   Suppose $H$ is a 
homogeneous polynomial of degree  $e$,  and suppose that $r_P(H) \geq 
r_P(C) + r_P(D) - 1$  for all points $P$ in and infinitely near to $C$ and $D$.  
Then there is an equation  
\[
H  =  A\cdot F + B\cdot G \,,
\]
where $A$  and  $B$  are homogeneous polynomials of degrees $e-c$  and  
$e-d$.  
\end{theorem}
\begin{proof}  By Corollary \ref{C2.6}, at every point $P$ in the plane, a local 
equation for $H$ in  $\mathcal{O}_P(\mathbb{P}^2)$  is in the ideal generated 
by local equations for $F$ and $G$. The fact that this is true locally if and only if 
it is true globally, so that there is an identity as shown, is proved in \cite{F}, 
\S5.5.
\end{proof}

It follows from what we have done here that most of the results proved in 
\cite{F} for plane curves with only ordinary singularities extend without 
essential change to curves with arbitrary singularities.  For example, suppose 
$C$ is an irreducible plane curve, defined by a homogeneous polynomial  
$F(X,Y,Z)$,  and  $\mathcal{X} = \widetilde{C} \to C$  is its nonsingular model 
(constructed by a 
succession of blowups over singular points of  $C$),  the {\em 
adjoint divisor}  $\Delta$  is the sum  $\sum d_QQ$,  where $d_Q$  is the 
order of vanishing at  $Q$  in  $\widetilde{C}$  of the conductor ideal at the 
image of  $Q$  in $C$.   The {\em genus}  $g_{\mathcal{X}}$  of 
$\mathcal{X}$ can be defined to be  
$(n-1)(n-2)/2 - \delta$,  where  $\delta = \sum\delta_P = 
\frac{1}{2}\operatorname{deg}(\Delta)$.  We choose coordinates so that the 
line $Z = 0$  intersects $C$ only at nonsingular points.  
	
The {\em divisor}  $\operatorname{div}(\omega)$  of a differential  $\omega$  
of  $R = R(\mathcal{X}) = R(C)$  over  $k$  is  
$\sum\operatorname{ord}_Q(\omega)Q$,  the sum over the (finite) set of  
$Q$  in $\mathcal{X}$ 
at which the order of  $\omega$  is not zero.  If  $f = F(x,y,1)$,  it 
follows from Proposition \ref{P3.1} that the order of  $dx/f_y$  at each point of 
the affine plane  $\mathbb{A}^2$  is  $-d_Q$.  One can calculate the order at 
the points of  $Z = 0$  by changing coordinates from the given copy of  
$\mathbb{A}^2$  to the other two copies of  $\mathbb{A}^2$.  One finds that  
\[
\operatorname{div}(dx/f_y)  =  -\Delta + (n-3)\operatorname{div}(Z) \, .
\]
A homogeneous polynomial $G(X,Y,Z)$  is {\em adjoint} to $C$ if the divisor  
$\operatorname{div}(G)$  cut out on $\mathcal{X}$ 
by $G$  contains  $\Delta$,  i.e.,  
$\operatorname{div}(G) = 
\Delta + A$,  for some effective divisor $A$.   If $G$  
is an adjoint to $C$ of degree  n-3,  it follows that  $\operatorname{div}(G) = 
\Delta + A$,  where  $A = \operatorname{div}(\omega)$  for some everywhere 
regular differential  $\omega$  on $\mathcal{X}$, 
namely $\omega = (G/Z^{n-3})\, dx/f_y$.  
Such adjoints exist whenever the 
genus is positive, since the condition for $G$  to be in the adjoint ideal at $P$ 
is defined by  $\delta_P$  linear equations, and the projective space of such 
forms has dimension  $(n-1)(n-2)/2$.  
	
The classical proof of the Riemann-Roch theorem, given in \cite{F}, Chap. 8, 
for curves with ordinary singularities, then applies without change for curves 
with arbitrary singularities.  This proof is based on Max Noether's 
Fundamentalsatz.  In particular, one sees that the adjoints of degree  $n-3$  cut 
out, besides the fixed component  $\Delta$,  the complete linear series of 
canonical divisors.  See \cite{ACGH}, Chap. I, App. A, for a modern discussion of 
adjoints and differentials for complex curves; there Gorenstein's theorem is 
deduced from the Riemann-Roch theorem.  Zariski (\cite{Z}, \S15) discusses adjoints 
in higher dimensions.  One can find a comparison with other notions of adjoints 
in \cite{GV1} and \cite{GV2}, and more about adjoints and conductors in 
\cite{AS}.

\section*{Appendix.  Resolution of singularities for planar curves}

We keep the notation of Sections 1 and 2.  We show that, for planar curves, the 
blowing up process must stop, by induction on the multiplicity.  Note that for 
one blowup,  
\[
\sum r_{P_i} \leq   \dim(A'/xA')  =  \dim(A/xA)  =  r_P  =  r\, , 
\]
by Remark \ref{R2.7}.  
It therefore suffices to show that it is impossible for there 
to be, at every stage, just one point on the proper transform of the curve over 
$P$,   with the same multiplicity  $r > 1$.   
	
We rule this out by a power series calculation.  We may assume that the leading 
term of $F$ is  $Y^r$.   Set  $F^{(1)} = F$,  and construct inductively a 
sequence of polynomials  $F^{(n)} = F^{(n)}(X,Y)$,  each of whose leading terms 
is  $Y^r$,  and a sequence of elements  $a_n$  in the ground field  $k$,  such 
that  $F^{(n-1)}(X,XY) = X^rF^{(n)}(X,Y-a_nX)$.\footnote{For example, if  
$F(X,Y) = Y^2 +2X^2Y + X^4 + X^7$,  then  $F(X,XY) = X^2((Y+X)^2 + X^5)$,  
so  $F^{(2)}(X,Y) = Y^2 + X^5$,  and  $a_2 = -1$;  then  $F^{(3)} = Y^2 + X^3$,  
with  $a_3 = 0$.}   It follows by induction that for all  $n \geq 2$,  with  
$\varphi_n(X) =  \sum_{i=2}^n a_iX^i$,
\[
F(X,X^{n-1}Y+\varphi_n(X))  =  X^{r(n-1)}F^{(n)}(X,Y)\, .
\]
Setting  $\varphi(X) =  \sum_{i=2}^{\infty} a_iX^i$  in  $k[[X]]$,  
we see that  
$F(X,\varphi(X)) = 0$, and so  $Y - \varphi(X)$  divides  $F(X,Y)$.  
If  $\varphi(X)$  
is a polynomial, this contradicts 
the irreducibility of $F$, so an infinite number 
of its coefficients must be nonzero.  We claim that  
$F(X,Y) = (Y - \varphi(X))^r$.  
If not, write  $F(X,Y) = (Y - \varphi(X))^s\cdot 
G(X,Y)$,  for some $G(X,Y)$  in  $k[[X]][Y]$  with  $G(X,\varphi(X)) \neq 0$,  
and  $s < r$.  From the displayed equation we see that  
\[
\begin{aligned}
(X^{n-1}Y+\varphi_n(X)-\varphi(X))^s \cdot 
&G(X,X^{n-1}Y+\varphi_n(X))  \, = \notag  \\
&X^{r(n-1)}F^{(n)}(X,Y) \, . \notag  
\end{aligned}
\]
If  $a_{n+1} \neq 0$,  setting  $Y = 0$  in this equation and computing the 
order of vanishing with respect to $X$,  one sees that  $G(X,\varphi_n(X))$  is 
divisible by  $X^{n(r-s) - s}$.   Since there are arbitrarily large  $n$  with  
$a_{n+1} \neq 0$, this shows that  $G(X,\varphi(X)) = 0$,  a contradiction. 
	
To complete the proof it remains to verify that if  $(Y - \varphi(X))^r$  is a 
polynomial, then  $\varphi(X)$  must be a polynomial.  This is clear in 
characteristic zero, so assume the characteristic is  $p$,  and write  $r = 
q\cdot u$,  with $q$  a power of  $p$  and $u$  relatively prime to  $p$.  Since 
the binomial coefficient $\binom{r}{q}$ is not zero modulo  $p$,  
$\varphi(X)^q$  must be in  $k[X]$,  and this implies that  $\varphi(X)$  is in  
$k[X]$.

This calculation shows that if  $\Lambda$  is the localization of  $k[X,Y]$  at its 
maximal ideal  $(X,Y)$,  and $F$ is any irreducible element in  $\Lambda$,  
then $F$ cannot be a power of an irreducible element in the completion  
$k[[X,Y]]$  of  $\Lambda$;  in other words, the completion of the ring  $A = 
\Lambda/(F)$  cannot have nilpotent elements.  
This illustrates the general fact that the 
integral closure of a one-dimensional Noetherian domain $A$  is a finitely 
generated  $A$-module if and only if its completion has no nilpotents; see 
\cite{N}, \S33, \cite{M1}, \S31, or \cite{M2}, \S33.   
	
An example from Nagata \cite{N}, Appendix, shows that this is not true for all 
two dimensional regular local rings  $\Lambda$.  To see such an example, let  
$\{a_{ij} \mid i, j \geq 0\}$  be a collection of indeterminates over 
$\mathbb{F}_p$.  Let  $K = \mathbb{F}_p({a_{ij}})$  be the field generated 
over  $\mathbb{F}_p$  by these indeterminates, and let  $\Lambda$  be the 
subring of  $K[[X,Y]]$  consisting of power series whose coefficients lie in some 
finite extension of  $K^p = \mathbb{F}_p({a_{ij}{}^p})$.  This  $\Lambda$  is a 
regular local ring, with maximal ideal generated by $X$ and $Y$,  and $F = 
\sum a_{ij}{}^pX^{pi}Y^{pj}$  is an element of  $\Lambda$  which is a  
$p^{\text{th}}$  power in the completion of  $\Lambda$,  but $F$ is not a  
$p^{\text{th}}$  power in  $\Lambda$.  One can verify directly that the blowing 
up process on this $F$ continues indefinitely.


\begin{thebibliography}{22}
%
\addcontentsline{toc}{section}{References}

\bibitem{AS} S.S.~Abhyankar and A.M.~Sathaye, {\em Geometric Theory of 
Algebraic Space Curves}, Springer Lecture Notes {\bf 423}, 1974

\bibitem{ACGH} E.~Arbarello, M.~Cornalba, P.A.~Griffiths, and J.~Harris, 
Geometry of Algebraic Curves, Vol. I. Springer-Verlag, 1985

\bibitem{D}  P.~Deligne, Intersections sur les surfaces r\'{e}guli\`{e}re, 
Expos\'{e} X in {\em S\'{e}minaire de G\'{e}om\'{e}trie Alg\'{e}brique du 
Bois-Marie},
by P.~Deligne and N.~Katz, Springer Lecture Notes {\bf 340}, 1973, 
pages 1--38

\bibitem{E}  D.~Eisenbud, {\em Commutative Algebra with a View Toward 
Algebraic Geometry}, Springer-Verlag, 1995

\bibitem{F}  W.~Fulton,  {\em Algebraic Curves, An Introduction to Algebraic 
Geometry}, W. A. Benjamin, Inc., 1969; new edition in preparation

\bibitem{GV1} S.~Greco and P.~Valabrega, On the theory of adjoints, in 
{\em Algebraic Geometry, Proceedings, Copenhagen 1978}, Springer Lecture 
notes {\bf 732} (1979), 98--123

\bibitem{GV2}  S.~Greco and P.~Valabrega, On the theory of adjoints II, 
Rend. Circ. Mat. Palermo {\bf 31} (1982), 5--15

\bibitem{G}  D.~Gorenstein, An arithmetic theory of adjoint plane curves, 
Trans. Amer. Math. Soc. {\bf 72}, (1952), 414--436

\bibitem{K}  E.~Kunz, {\em Ebene algebraische Kurven}, Der Regensburger 
Trichter {\bf 23}, Fakult\"{a}t f\"{u}r Mathematik der Universit\"{a}t 
Regensburg, 1991

\bibitem{LJ}  M.~Lejeune-Jalabert, Le Th\'{e}or\`{e}me ``AF + BG'' de Max 
Noether, in {\em S\'{e}minaire sur les Singularit\'{e}s}, L\^{e} D\u{u}ng 
Tr\'{a}ng, Publications Math\'{e}matiques de l'Universit\'{e} Paris VII, 1980, 
pages 97--138

\bibitem{M1} H.~Matsumura, {\em Commutative Algebra}, Second edition, 
Benjamin/Cummings, 1980.

\bibitem{M2}  H.~Matsumura, {\em Commutative Ring Theory}, Cambridge 
University Press, 1986

\bibitem{N}  M.~Nagata, Local Rings, Interscience Publishers, 1962

\bibitem{Sa}  P.~Samuel, Singularit\'{e}s des vari\'{e}t\'{e}s alg\'{e}briques, 
Bull. Soc. Math. France {\bf 79}, (1951), 121--129

\bibitem{Se}  J.-P.~Serre, {\em Algebraic Groups and Class Fields}, 
Springer-Verlag, 1988

\bibitem{Sz}  L.~Szpiro,  {\em Lectures on equations defining space curves}, 
Notes by N. Mohan Kumar, Tata Institute of Fundamental Research, Bombay, 
Springer-Verlag, 1979

\bibitem{Z} O.~Zariski, {\em An Introduction to the Theory of Algebraic 
Surfaces}, Springer Lecture Notes {\bf 83}, 1969

\end{thebibliography}
\end{document}